\documentclass{article}

\usepackage{microtype}
\usepackage{graphicx}
\usepackage{subfigure}
\usepackage{booktabs} 
\usepackage{amsmath}
\usepackage{mathtools}
\usepackage{subfigure}
\usepackage{caption}
\usepackage{algorithm}
\usepackage{algorithmic}
\usepackage{hyperref}
\newcommand{\maxm}{\mathop{\mathrm{maximize}}}

\usepackage[archival]{icml2023}

\usepackage{amsmath}
\usepackage{amssymb}
\usepackage{mathtools}
\usepackage{amsthm}

\usepackage[capitalize,noabbrev]{cleveref}

\theoremstyle{plain}
\newtheorem{theorem}{Theorem}[section]

\theoremstyle{definition}

\theoremstyle{remark}

\newtheorem{example}[theorem]{Example}
\usepackage[textsize=tiny]{todonotes}

\icmltitlerunning{Submission and Formatting Instructions for ICML 2023}

\begin{document}

\twocolumn[
\icmltitle{ReLU Neural Networks, Polyhedral Decompositions, and Persistent Homology}




\begin{icmlauthorlist}
\icmlauthor{Yajing Liu}{yyy}
\icmlauthor{Christina Cole}{yyy}
\icmlauthor{Christopher Peterson}{yyy}
\icmlauthor{Michael Kirby}{yyy}
\end{icmlauthorlist}

\icmlaffiliation{yyy}{Department of Mathematics, Colorado State University, Fort Collins, USA}

\icmlcorrespondingauthor{Yajing Liu}{yajing.liu@colostate.edu}

\icmlkeywords{Machine Learning, ICML}

\vskip 0.3in
]



\printAffiliationsAndNotice{}  

\begin{abstract}
A ReLU neural network leads to a finite polyhedral decomposition of input space and a corresponding finite dual graph. We show that while this dual graph is a coarse quantization of input space, it is sufficiently robust that it can be combined with persistent homology to detect homological signals of manifolds in the input space from  samples. This property holds for a variety of networks trained for a wide range of purposes that have nothing to do with this topological application. We found this feature to be surprising and interesting; we hope it will also be useful.
\end{abstract}

\section{Introduction}
\label{intro}
The rectified linear unit (ReLU)  function is the default activation function for many well known, deep, feed forward  neural networks (AlexNet, ResNet, Inception, SqueezeNet, etc). These networks frequently use some convolutional layers to decrease the number of parameters being trained and often utilize skips to decrease issues from overfitting.
The importance of ReLU feed forward neural networks (FFNNs) has motivated an extensive investigation of their properties from a variety of aspects to demystify their ``black-box nature".
\textcolor{black}{In this paper, we start from the observation that, regardless of their use of skips and/or convolutions, a ReLU FFNN $F:\mathbb R^m \rightarrow \mathbb R^n$ decomposes the input space ($\mathbb R^m$) into convex polyhedra and assigns to each polyhedron a unique binary vector that encodes the ReLU activation pattern of all nodes in the ReLU layers of the network.}
\textcolor{black}{More precisely, the neural network assigns a binary vector to each input point, assigning identical binary vectors to input points lying in the same polyhedra.} 
\textcolor{black}{Despite this granularity imposed on input space, we found that the homology of manifolds embedded in the input space of the network can still be detected via persistent homology applied solely to distance measures built from the binary vectors of the points in the sampling.} After explaining the connection between polyhedra and binary vectors in a ReLU FFNN, we describe methods for finding the binary vectors associated with the nearby neighbors of a given polyhedron in the decomposition. For small networks, one can (in principle) use these methods to determine all the polyhedra in the decomposition and determine their proximity within the decomposition. Toward this goal, we exploit the following observation: {\it two polyhedra share a facet if and only if their associated binary vectors differ in exactly one bit}. Due to this property, we propose the Hamming distance between binary vectors as a proxy for proximity between polyhedra in the decomposition. More precisely, if $\mathcal{X}$ is a collection of data points in the input space, each $x\in \mathcal{X}$ lies in one of the polyhedra and can be assigned the binary vector, $s(x)$, that labels the polyhedron. Using the Hamming metric, we build a network driven distance matrix between pairs of binary vectors associated to data points in $\mathcal{X}$. Using a software package such as Ripser (see {\tt https://live.ripser.org}), we illustrate how the coarse distance measure corresponding to a polyhedral decomposition of input space can be used to extract homological information about a manifold lying in the input space from sample points on the manifold. It is worth noting that the Hamming distance between polyhedral bit vectors is an approximation for the smallest number of polyhedral steps between two polyhedra (i.e. for the shortest path between vertices on the dual graph of the polyhedral decomposition). \cite{Jamil_2023_CVPR} has demonstrated the efficacy of the Hamming distance in illuminating the mechanisms underlying neural network functionality.

The paper is organized as follows: 
Section \ref{sec:not-and-def} introduces the basic background on ReLU FFNNs, the polyhedral decomposition associated to a network, the binary vectors assigned to each polyhedron, the linear model that determines each polyhedron, the linear program for filtering out redundant  inequalities in the linear model, and the affine linear function that is attached to each polyhedron. Section \ref{sec:alg-bit-search} outlines  algorithms for determining the binary vectors that occur for polyhedra in the entire input space and for polyhedra that occur within a bounded region. Section \ref{sec:pers-hom} describes how the polyhedral decomposition, realized through binary vectors, can be combined with persistent homology to uncover topological features in data. Section \ref{sec:concl} concludes the work and sketches directions for future research.

\section{Neural Networks and Polyhedral Decompositions}\label{sec:not-and-def}
This section gives the basic definitions of a ReLU neural network and describes the connection to polyhedral decompositions and binary vectors.
\subsection{Notation for Neural Networks}
Consider an $(L+1)-$layer ReLU feedforward neural network: 
\begin{multline}
    \label{eq:nn}
\mathbb{R}^m\xrightarrow[\text{ReLU}]{(W_1, b_1)}\mathbb{R}^{h_1}\xrightarrow[\text{ReLU}]{(W_2, b_2)}\mathbb{R}^{h_2}\rightarrow\ldots\rightarrow\mathbb{R}^{h_{L-1}}\\
\xrightarrow[\text{ReLU}]{(W_{L}, b_{L})}\mathbb{R}^{h_L}
\xrightarrow[]{(W_{{L+1}}, b_{{L+1}})}\mathbb{R}^{n}.
\end{multline}

In this model, $\mathbb{R}^m$ is the input space, $\mathbb{R}^n$ is the output space, and $h_i$ corresponds to the number of nodes at layer $i$. Layer 0 corresponds to the input space and layer $L+1$ corresponds to the output space (so $h_0 = m$ and $h_{L+1}=n$). We let $W_i\in \mathbb{R}^{h_i\times h_{i-1}}$ and  $b_i\in\mathbb{R}^{h_i}$ denote the weight matrix and bias vector of layer $i$, {respectively}. The activation function{s} for the hidden layers (layers $1, \dots, L$)  are assumed to be ReLU functions (applied coordinate-wise) while 
the map to the last layer (the output layer) is assumed to be affine linear (without a ReLU function being applied to the image). Recall that the ReLU function is the map $ReLU: \mathbb R \rightarrow \mathbb R$ given by
\begin{equation}
      ReLU(a)=  
    \begin{cases}
      a & \text{if}\  a>0 \\
      0 & \text{if}\  a\leq 0.
    \end{cases} 
\end{equation}

The ReLU map is a map on real numbers that is piecewise linear and continuous.  It can be naturally extended to a piecewise linear continuous map on vector spaces (which we also denote as ReLU). More precisely, we define $ReLU: \mathbb R^{h_i} \rightarrow \mathbb R^{h_i}$, by applying the ReLU function to each coordinate of $x\in \mathbb R^{h_i}$.
Let $w_{i,j}$ denote the $j^{th}$ row of $W_i$ and let $b_{i,j}$ denote the $j^{th}$ entry of $b_i$. 
Given an input data point $x\in\mathbb{R}^m$, we denote the output of $x$ in layer $i$ as $F_i(x)$. Thus, with this notation we have $F_i(x)\in\mathbb{R}^{h_i}$,  $F_0(x) = x$, and 
\begin{equation}
    \begin{aligned}
\label{eq:outputderi}
  F_i(x)&=\text{ReLU}(W_iF_{i-1}(x)+b_{i})\\
    &= \begin{bmatrix} \max\{0, w_{i,1}F_{i-1}(x)+b_{i,1}\}\\
    \vdots\\
    \max\{0, w_{i, h_i}F_{i-1}(x)+b_{i, h_i}\}
    \end{bmatrix}.   
\end{aligned}
\end{equation}

\subsection{Definitions of Binary Vectors}
\label{subsec:bitvec}
Consider model (\ref{eq:nn}). Given an input data point $x\in\mathbb{R}^m$, for each hidden layer $i$ (so $1\leq i\leq L)$, we introduce a binary (bit) vector \[s_i(x)=[s_{i,1}(x) \ \ldots\  s_{i,h_i}(x)]^{\top}\in\mathbb{R}^{h_i},\] where $s_{i, j}(x)$ (with $1\leq j\leq h_i$) is defined as follows:
\begin{equation}
\label{eq:bitvector}
        s_{i, j}(x)=  
    \begin{cases}
      1 & \text{if}\  w_{i,j}F_{i-1}(x)+b_{i,j}>0 \\
      0 & \text{if}\  w_{i,j}F_{i-1}(x)+b_{i,j}\leq 0.
    \end{cases} 
\end{equation}
Thus, for each point $x\in \mathbb R^m$, we have a sequence of binary vectors $s_1(x), s_2(x), \dots, s_L(x)$.
We can stack the binary vectors associated to $x$  to make a long column vector 
\begin{equation}
\label{eq:bitvec}
    s(x)= [s_1^{\top}(x)\  \ldots \ s_{L}^{\top}(x)]^{\top}\in\mathbb{R}^{h},
\end{equation}
where $h = \sum_{i=1}^L h_i$ is the total number of nodes in the hidden layers. We call $s(x)$  the binary vector of $x$.

Different points from the input space $\mathbb{R}^m$ can have the same binary vector. We next show that the set of points that have the same binary vector, $\{x' :  s(x')=s(x), x' \in \mathbb{R}^m\}$, form a convex polyhedron in $\mathbb{R}^m$.

\subsection{Linear Model for Binary Vectors}
Let $s_1, s_2, \dots, s_L$ denote a given sequence of binary vectors for model (\ref{eq:nn}). For each layer $i$ $(1\leq i\leq L)$, we describe inequality constraints that must be satisfied for an input data point to have the same binary vector $s_i$.
To describe the inequalities in a consistent manner, we introduce a sign vector of $1's$ and $-1's$ for each hidden layer.
For layer $i$, define $s_i'=[s'_{i,1} \ \ldots \ s'_{i,h_i}]^{\top}$ with
        \begin{equation}
        s'_{i, j}=  
    \begin{cases}
      1 & \text{if}\  s_{i,j}=0\\
      -1 & \text{if}\ s_{i,j}=1.
    \end{cases}       
\end{equation}
\noindent
In layer 1, because $F_1(x)=\text{ReLU}(W_1x+b_1)$, any data point $x$ that has the bit vector $s_1$  satisfies the following linear inequality:
\begin{equation}
    \label{ineq:layer1}
    \text{diag}(s_1')({W}_{1}x+b_1) \leq 0 
\end{equation}
{where diag$(v)$ is a square diagonal matrix with the elements of vector $v$ on the main diagonal. The inequality (\ref{ineq:layer1}) is established by a fundamental principle: the affine output generated by the first hidden layer, in response to input $x$, must satisfy a greater-than-zero condition for nodes where the corresponding binary vector is set to 1, and a less-than-or-equal-to-zero condition for nodes where the corresponding binary vector is set to 0.}

 Suppose that $x$ has the bit vector sequence $s_1, s_2, \dots, s_L$.
Let $\hat{W}_j=W_j\text{diag}(s_{j-1})\hat{W}_{j-1}$
and $\hat{b}_j = W_j\text{diag}(s_{j-1})\hat{b}_{j-1}+b_j$ for $2\leq j\leq L$ with $\hat{W}_1 = W_1, \hat{b}_1 = b_1$. By model (\ref{eq:nn}),  the following equation holds {for $1\leq j\leq L$:
\begin{multline*}
           F_j(x) =\text{ReLU}(W_jF_{j-1}(x)+b_j) = \text{diag}(s_j)(\hat{W}_jx+ \hat{b}_j)
        \end{multline*}
where $F_0(x)=x$.}

{More generally, any data point $x$ that has $s_j$ as its bit vector for layer $j$   should satisfy the following linear inequalities:
\begin{equation}
    \label{ineq:layerj}
    \text{diag}(s_j')\hat{{W}}_{j}x \leq \text{diag}(s_j')(-\hat{b}_j).
\end{equation}
\noindent
Let $A_j = \text{diag}(s_j')\hat{W}_{j}$ and $c_j=\text{diag}(s_j')(-\hat{b}_j)$ for $1\leq j\leq L$. Combining (\ref{ineq:layer1}) and (\ref{ineq:layerj}), we have 
\begin{equation}
    \label{ineq:polytopecond}
    Ax\leq c,
\end{equation}
where {\[A=[A_1^{\top}\ A_2^{\top} \ldots \ A_L^{\top}]^{\top}\ \text{and}  \ c=[c_1^{\top}\ c_2^{\top} \ldots \ c_L^{\top}]^{\top}\]}
with $A_i\in \mathbb{R}^{h_i\times m}$ and $c_i\in\mathbb{R}^{h_i}$. Note that the set $P=\{x\in\mathbb{R}^m: Ax\leq c\}$ is a convex polyhedron (potentially empty, potentially unbounded). Considering the totality of points in the input space, the weights of a ReLU neural network lead to a decomposition of the input space into a collection of bounded and unbounded polyhedra and each polyhedron has a corresponding bit vector. Note that a random bit vector may or may not correspond to a non-empty polyhedra.

It is not difficult to see that for a polyhedron $P$ of full dimension (the same dimension as the ambient space), there is a unique smallest subset of inequalities, which we denote by $(A',c')$, that one can obtain from $(A,c)$ that leaves the polyhedron $P$ unchanged. Thus, $A'$ is built from a subset of rows of $A$ and $c'$ is the corresponding subset of rows of $c$ such that the set of points that satisfy $Ax\leq c$ is the same as the set of points that satisfy $A'x\leq c'$. 

Write \[A=\begin{bmatrix} a_1\ a_2\ \ldots \ a_h \end{bmatrix}^{\top}\  \text{and}\quad  c = \begin{bmatrix}
    c^1\ c^2\ \ldots \ c^h
\end{bmatrix}^{\top}\] with $a_i\in\mathbb{R}^m$ and $c^i\in\mathbb{R}$. 
To determine if {$a_ix\leq c^i$} is a redundant constraint, first define \[\Tilde{A}=\begin{bmatrix} a_1\ a_2\ \ldots \ a_{i-1}\ a_{i+1} \ldots\ a_h \end{bmatrix}^{\top} \] \text{and} \[\Tilde{c}=\begin{bmatrix}
    c^1\  c^2\  \ldots \ c^{i-1}\ c^{i+1}\ \ldots\ c^h
\end{bmatrix}^{\top}.\] Next, consider the following linear program
\begin{align}\label{prob:LP}
\begin{array}{l}
\maxm\limits_{x} \ \ a_i^{\top}x \\
\text{s.\ t.} \ \ \Tilde{A}x\leq \Tilde{c}.
\end{array}
\end{align}
If the optimal objective value of (\ref{prob:LP}) is less than or equal to $c^i$, then the $i$th linear inequality is redundant, and we can remove $a_i^{\top}$ and $c^i$ from $A$ and $c$, respectively.  We determine $(A', c')$ by iterating this process to remove all redundant constraints.

Given a bit vector $s$, the $i$-th entry of $s$ is called active if the $i$th row of $A$ is in $A'$, and inactive otherwise. For any input $x$ in the polyhedron determined by the bit vector $s$, the output of $x$ is determined by the single affine map: $G(x)=W_{L+1}\text{diag}(s_L)\hat{W}_Lx+ W_{L+1}\text{diag}(s_L)\hat{b}_L+b_{L+1}$.

We summarize some of the points of this section, together with implications (many of which we leave to the reader) in the bullet points that follow. We emphasize that the polyhedral decomposition of input space has an associated dual graph. This graph has vertices corresponding to polyhedra and edges corresponding to polyhedra that share an $m-1$ dimensional face. It is not hard to show that this dual graph is bipartite (refer to Appendix \ref{appendix:2Dpoly}). We will be using this graph, implicitly, to detect homological signals via persistent homology. For the interested reader, we point out the following references related to topological features and combinatorial features of the type of neural networks we are considering \cite{grigsby2022local, masden2022algorithmic}. We would also like to acknowledge the quickly-growing library of work in the area of polyhedral theory, a survey of which can be found in \cite{huchette_when_2023}.

The following is a summary of the key ideas presented either explicitly or implicitly in this section:
\begin{itemize}
    \item A ReLU FFNN determines a continuous piecewise linear map $F:\mathbb R^m \rightarrow \mathbb R^n$. \vspace{-0.70em}
    \item The neural network leads to a decomposition of $\mathbb R^m$ into a collection of bounded and unbounded convex polyhedra.\vspace{-0.70em}
    \item A binary vector, supported on the hidden nodes of the network, can be attached to each point in the domain of the network. 
    The value of the binary vector associated to a given input
    $x \in \mathbb{R}^m$,
    at a given node, is $1$ if ReLU was not applied at the node and $0$ if ReLU was applied at the node (the ReLU activation pattern). \vspace{-0.70em}
    \item If $h$ denotes the number of hidden nodes in the network, then there are, a priori, $2^h$ possible binary vectors and thus $2^h$ possible polyhedra in the decomposition of input space. In reality, the number of realized convex polyhedra in the domain is much smaller. \vspace{-0.70em}
    \item If two points lie in the interior of the same convex polyhedron then they determine the same binary vector.\vspace{-0.70em}
    \item If two points are in the interior of distinct convex polyhedra then they determine distinct binary vectors. \vspace{-0.70em}
    \item From the ReLU activation pattern for points in a convex polyhedron, one can determine an affine linear equation that represents the behavior of the neural network on points in the polyhedron. If the polyhedron is denoted by $P$, then for points in $P$, the function $F:\mathbb R^m \rightarrow \mathbb R^n$ can be expressed as $F(x)=A_Px+b_P$ for some matrix $A_P$ and some vector $b_P$. \vspace{-0.70em}  
    \item Input values that determine a value of $0$ on some node, before applying ReLU, are precisely the input values that lie on the boundary of distinct convex polyhedra. This is due to the fact that applying ReLU to $0$ has the same effect as not applying ReLU to $0$. \vspace{-0.70em}
    \item If two convex polyhedra, $P_1,P_2$, share an $(m-1)$-dimensional face, then the binary vectors associated to each of the polyhedra differ in one bit. Furthermore, the affine linear functions $A_{P_1}+b_{P_1}$ and $A_{P_2}+b_{P_2}$ agree on this $(m-1)$-dimensional face. \vspace{-0.70em}
    \item Any polyhedral decomposition of $\mathbb R^m$ has a natural dual graph with vertices corresponding to $m$-dimensional polyhedra and edges corresponding to polyhedra sharing an $m-1$ dimensional facet. \vspace{-0.70em}
    \item The Hamming distance between binary vectors (that represent two polyhedra) can be used as an approximation for the smallest number of polyhedral steps between the two polyhedra (i.e. the length of a minimal geodesic on the dual graph).
    \end{itemize}

\section{Algorithms for Bit Vector Search and Examples}\label{sec:alg-bit-search}

The number of  polyhedra in the input space or within a bounded region of the input space provides a measure of the network's expressivity and complexity. Upper
and lower bounds on the maximal number of polyhedra obtainable from a given ReLU FFNN architecture can be found in
 \cite{pascanu_number_2013}, \cite{montufar_number_2014}, \cite{raghu_expressive_2017}, \cite{arora2018understanding}, and \cite{serra_bounding_2018}.  Several algorithms (\cite{xiang_reachable_2018}, \cite{yang_reachability_2020}, and \cite{xu_traversing_2021}) have been developed to compute the exact polyhedra decomposition of the input space through \textcolor{black}{layer-by-layer} linear inequality solving. Larger decomposition examples can be computed using a method developed by 
\cite{vincent_reachable_2021}, which
enumerates all polyhedra in the input space as follows:
\begin{itemize}
    \item Start with a random point $x\in \mathbb R^m$ and determine its bit vector $s(x)$. This bit vector labels a polyhedron $P$. \vspace{-0.70em}
    \item Find the active bits in $s(x)$ for the polyhedron $P$. \vspace{-0.70em}
    \item Each active bit corresponds to a neighboring polyhedron. Each neighboring polyhedron has a bit vector that can be obtained by ``flipping'' one of the active bits for $P$. Thus, one can find all of the neighboring polyhedra for $P$, in terms of their binary vectors. \vspace{-0.70em}
    \item Repeat the process to find the neighbors of each of these newly identified polyhedra.\vspace{-0.70em}
    \item The number of active bits for a polyhedron is equal to the number of nearest neighbors of the polyhedron. The previous steps continue until we have a list of polyhedra $\mathcal{P}$ 
    that satisfies the property that for each $P\in \mathcal{P}$ 
    , the set of nearest neighbors of $P$ is itself a subset of $\mathcal{P}$. 
\end{itemize}
This process leads to a set $\mathcal{P}$ 
of convex polyhedra that decompose input space. Such decompositions have been used to define network-imposed distance metrics for model inputs \cite{balestiero_spline_2018}. We coarsen these metrics and show they can still be used to detect homological signals. Equating neural networks of various types (convolutional neural networks, residual networks, skip-connected networks, recurrent neural networks) to max-affine spline operators has been carried out by \cite{balestiero_spline_2018}. The paper \cite{sattelberg_locally_2020} investigated the behavior of the ReLU FFNNs based on the structure of the decomposition, together with the affine map attached to each polyhedron. 

While there are multiple methods that can be used to determine the polyhedral decomposition, and their associated binary vectors, imposed by the weights of a ReLU neural network, all of these methods are woefully inadequate for the large and deep networks that are commonplace today. This is not a fault of the algorithms, many modern networks contain many millions of hidden nodes and have input spaces with dimension well beyond 100,000. There are just too many polyhedra represented by such networks.  However, for small networks, any of the methods will work. We present two straightforward methods. The first is a brute-force method formulated as a Linear programming problem. The second, which is the traversal-and-pruning method outlined above, computationally improves upon the first. We demonstrate our methods with examples. Pseudocode for these two methods are provided in Appendix~\ref{appendix:algorithms}.

\subsection{Algorithms}
 A mixed-integer linear program to count the number of polyhedra in a bounded region can be found in \cite{serra_bounding_2018}. {We present a linear program that can also count polyhedra, determine if they share a facet, and whose implementational simplicity is useful. Our proposed method yields not only the count of polyhedra in the entire input space $\mathbb{R}^m$ but also their respective binary vectors}. To count the number of polyhedra in a bounded region, the linear inequality constraints in (\ref{eq:LP}) need to be augmented to include bounds for variables. For example, when the input is 2-dimensional and the bounded region is determined by $a_1\leq x_1\leq b_1, \ a_2\leq x_2 \leq b_2$, the constraints become $\Tilde{A}^jx\leq \Tilde{c}^j$ where $\Tilde{A}^j$ is the concatenation of $A^j$ and $B$ while $\Tilde{c}^j$ is the concatenation of $c^j$ and $d$, with 
$$B=\begin{bmatrix}
    1 &-1 &0 &0 \\
   0 &0 &1 &-1
\end{bmatrix}^T \ \  {\rm and} \ \ d =\begin{bmatrix}
    b_1\ -a_1\ b_2\ -a_2
\end{bmatrix}^T.$$ 

\noindent
Algorithm~\ref{alg:bruteforce} is fast for small neural network structure. However,
\cite{hanin_deep_2019} proved that the number of bit vectors that correspond to polyhedra is approximately $h^m$ which is much smaller than $2^h$ for large $h$, highlighting the computational savings that can be procured by only traversing active bit vectors instead of all possibilities as is done by the brute-force method. These savings are what Algorithm \ref{alg:traversal_pruning} has to offer.

Algorithm 1 is easy to implement but slow.
The idea for Algorithm 2 is based on the fact that each polyhedron corresponding to $Ax\leq c$ is determined by  $A'x\leq c'$ and two adjacent polyhedra that share one facet differ by one active bit. This active bit corresponds to the hyperplane containing this facet. The key step for Algorithm 2 is to find the active bits for a given bit vector using the method mentioned in Section~\ref{subsec:bitvec}. Once identified, these active bits can then be used to find all neighboring polyhedra, starting a scheme that ripples through the desired domain until all polyhedra are found. The idea behind this algorithm also motivated the reachable polyhedral marching algorithm in \cite{vincent_reachable_2021}. We make a small modification to enumerate the polyhedra in a bounded region rather than in the entire, unbounded domain. 
For a bounded region, the bit vector can be expanded by checking whether the corresponding polyhedron hits the boundary or not. For example, when  $m =2$ and the bounded region is defined by $a_1\leq x_1\leq b_1, a_2\leq x_2 \leq b_2$, the bit vector $s$ can be expanded by 1 more bit with 1's or 0's by checking if its $Ax\leq c$ {solutions hit} any {of the domain boundaries} defined by $x_1= a_1$, $x_1= b_1$, $x_2= a_2$, or $x_2= b_2$ (1 if true, 0 otherwise). For the bit vectors that have $1$ for their last bit, we update their labels with 1 after they are added to $\mathcal{P}$. The reason is that only if the bit vector of the initial point does not hit the region boundaries, the bit vectors that hit the region boundaries will be derived by flipping the active bits of their neighbors, which means that all its adjacent polyhedra in the bounded region are already in $\mathcal{P}$.
\subsection{Examples}
In this section, we visualize the polyhedral decomposition, determined by a ReLU FFNN, on a bounded region for two basic models:
\begin{equation}
\small
    \label{example:model1}
    \mathbb{R}^2\xrightarrow[\text{ReLU}]{(W_1, b_1)}\mathbb{R}^{3}\xrightarrow[\text{ReLU}]{(W_2, b_2)}\mathbb{R}^3\xrightarrow[]{(W_3, b_3)}\mathbb{R}
\end{equation}
 \begin{equation}
 \small
     \label{example:model2}
     \mathbb{R}^3\xrightarrow[\text{ReLU}]{(W_1, b_1)}\mathbb{R}^{10}\xrightarrow[\text{ReLU}]{(W_2, b_2)}\mathbb{R}^{10}\xrightarrow[\text{ReLU}]{(W_3, b_3)}\mathbb{R}^{10}\xrightarrow[]{(W_4, b_4)}\mathbb{R}
 \end{equation}
We first applied model (\ref{example:model1})  to fit $f_1(x_1, x_2) =  x_1^2+x_2^2-2/3$ using 10000 points uniformly sampled in $[-1, 1]^2$ and model~(\ref{example:model2}) to fit $f_2(x_1, x_2, x_3) =  (x_1-1)^2+2x_2^2+x_3^2+1$ using 125000 points uniformly sampled in $[-1, 1]^3$. We used TensorFlow to train the above models with batch size of 50 and early stopping criterion based on the convergence of validation loss.

  We used Algorithm 1 to enumerate all binary vectors in $\mathbb R^2$ and in $[-1, 1]^2$. We used Algorithm 2 to enumerate all binary vectors in $\mathbb R^3$ and in $[-1, 1]^3$. We plotted the polyhedra by finding their extremal vertices. 

\begin{figure}[ht]
  \centering
  \subfigure[2-D polyhedral decomposition\label{fig:image-a-visualization}]{%
      \includegraphics[width=0.45\linewidth]
      {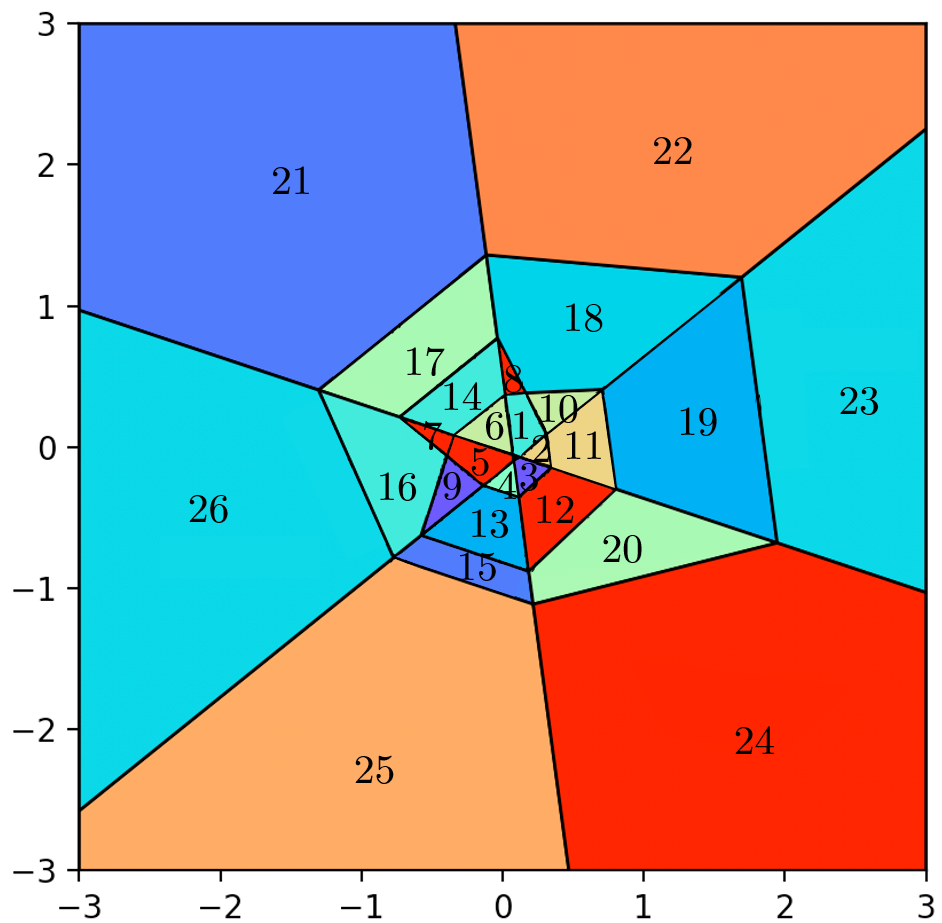}}
  \qquad
  \subfigure[3-D polyhedral decomposition \label{fig:image-b-visualization}]{%
      \includegraphics[width=0.45\linewidth]{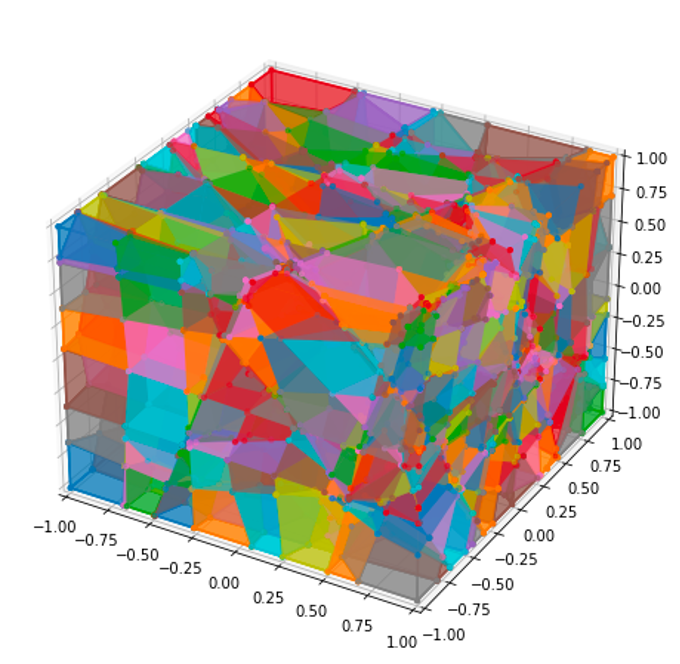}}
  \captionsetup{position=bottom}
  \caption{Visualizations of low-dimensional polyhedral decompositions. A detailed discussion on (a) can be found in Appendix \ref{appendix:2Dpoly} along with the exact binary vectors that label all of the regions of the decomposition. \label{fig:visualization}}
\end{figure}
\vspace{-0.3cm}
Figure~\ref{fig:visualization}(a) and Figure~\ref{fig:visualization}(b) provide visualizations of the polyhedral decomposition derived from model (\ref{example:model1}) and model (\ref{example:model2}), respectively. For Figure~\ref{fig:visualization}(a),
there are 25 polyhedra in the bounded region $[-1, 1]^2$  and 27 polyhedra in $\mathbb{R}^2$. Of these polyhedra, 21 are bounded and 6 are unbounded.  The binary vectors are superimposed onto the polyhedra. Note that the binary vectors associated to different polyhedra
differ in exactly one bit if and only if they share a facet. For Figure~\ref{fig:visualization}(b), there are 1858 polyhedra in the bounded region $[-1, 1]^3$ and  3331 polyhedra in $\mathbb{R}^3$.
\section{Persistent Homology and Polyhedral Decompositions}\label{sec:pers-hom}  {It has been well documented that the Euclidean distance between sampled points on a manifold in $\mathbb R^n$ can be employed to detect the topology of the manifold.  In this section, we provide a description of Vietoris-Rips persistent homology and illustrate how it can be effectively combined with the non-Euclidean distance measure, associated to the polyhedral decomposition, to also identify homological features.} Persistent homology has been a rapidly developing branch of topology largely due to its usefulness in data analysis and machine learning \cite{ghrist2008barcodes, zomorodian2004computing, carlsson2009topology, edelsbrunner2008persistent} (a collection of additional resources and videos can be found at {\tt https://www.aatrn.net}). Work linking persistent homology and neural networks has been appearing with increasing frequency, please see the following for a sampling of some recent works in this direction  \cite{rieck2018neural, zhao2020persistence, carriere2020perslay, birdal2021intrinsic}. The description of persistent homology given below is extremely condensed. The interested reader is encouraged to read the survey article \cite{ghrist2008barcodes} where additional details can be found.

\subsection{Persistent Homology}
Consider a metric space, $C$, with distance function $d: C\times C \rightarrow \mathbb R$. 
Let $x^{(1)}, x^{(2)}, \dots x^{(N)}$ be points in $C$. We can utilize $d$ to build an $N\times N$ matrix $D$ by setting $D_{i,j}=d(x^{(i)},x^{(j)})$. $D$ will be hollow (i.e. $D_{i,i} = 0$), symmetric, and non-negative. From $D$ we can build a family, $A(t)$, of $N\times N$ $\{0,1\}$-matrices parameterized by a real parameter $t$ using the rule
 \begin{equation}
        A(t)_{i,j}=  
    \begin{cases}
      0 & \text{if}\  D(i,j)>t \ \text{or\ if} \ i=j\\
      1 & \text{else}.
    \end{cases}       
\end{equation}
For each $t$, $A(t)$ can be viewed as an adjacency matrix of a graph $G(t)$. Let $CL(t)$ denote the clique complex of $G(t)$ \cite{hausmann_vietoris_1995}. By construction, $CL(t)$ is a simplicial complex and there is a natural inclusion map $$CL(t_1)\hookrightarrow CL(t_2)$$ whenever $t_1<t_2$. Let $\mathbb F$ be a field. A simplicial complex $S$ has an associated chain complex $S_{\bullet}$ of $\mathbb F$-vector spaces. The failure of the chain complex to be exact at location $i$ is measured by the $i^{th}$ homology $H_i( S_{\bullet},\mathbb F)$ (which itself is an $\mathbb F$-vector space). The inclusion map $CL(t_1)\hookrightarrow CL(t_2)$ induces a chain map $$CL(t_1)_{\bullet} \rightarrow CL(t_2)_{\bullet}.$$ Whenever you have a chain map $F:S_{\bullet} \rightarrow T_{\bullet}$ between chain complexes, $S_{\bullet}, T_{\bullet}$ you get associated linear maps between $H_i(S_{\bullet},\mathbb F)$ and $H_i(T_{\bullet},\mathbb F)$ for each $i$. Thus, the chain map $CL(t_1)_{\bullet} \rightarrow CL(t_2)_{\bullet}$ induces, for each $i$, a linear map $$H_i(CL(t_1)_{\bullet},\mathbb F) \rightarrow H_i(CL(t_2)_{\bullet},\mathbb F).$$ If we pick values $t_1<t_2< \dots <t_k$, then we can build a nested sequence of simplicial complexes 
$$\small CL(t_1)\subset CL(t_2) \subset \dots \subset CL(t_k)$$ which leads to \begin{equation} \small\label{PER}
    H_i(CL(t_1)_{\bullet},\mathbb F) \rightarrow H_i(CL(t_2)_{\bullet},\mathbb F) \rightarrow \dots \rightarrow H_i(CL(t_k)_{\bullet},\mathbb F),
    \end{equation} where the arrows denote linear maps.
    
    Each of the $H_i(CL(t_j)_{\bullet},\mathbb F)$ are finite dimensional $\mathbb F$-vector spaces and we can view (\ref{PER}) as a finitely generated graded $\mathbb F[x]$ module, where $\mathbb F[x]$ is the ring of polynomials. This module is frequently called the $i^{th}$ persistence module associated to the nested sequence of simplicial complexes. Since $\mathbb F[x]$ is a principal ideal domain, the $F[x]$ module has a decomposition into its invariant factors by the well known structure theorem for finitely generated graded modules over a principal ideal domain. Each invariant factor can be viewed as a homology class that has a birth time $t_b$ and a death time $t_d$ (possibly infinite) and this invariant factor can be represented as an interval $[t_b,t_d]$. This collection of intervals corresponds to the barcode representation of the invariant factors of a persistence module.

\subsection{Combining Persistence with a Polyhedral Decomposition}
First we provide two examples based on the polyhedral decomposition of model (\ref{example:model2}).

\begin{example} We consider the circle in $\mathbb R^3$ (the input space) with parameterization given by $(0,cos(t),sin(t))$. We sampled the circle at $20$ evenly spaced points $x^{(1)}, x^{(2)}, \dots, x^{(20)}$. Each of these points lie in one of the polyhedra from Figure~\ref{fig:visualization}(b). A picture of the polyhedra encountered by the $20$ sample points is found in Figure~\ref{fig:poly_sample}(a). We recorded the binary vectors for each of the encountered polyhedra. This gave a total of $20$ binary vectors but only $19$ were distinct. We labeled these distinct binary vectors {$s^{(1)}, s^{(2)}, \dots s^{(19)}$}. We built a $19\times 19$ matrix $E$ by setting $E_{i,j}$ equal to the Hamming distance between {$s^{(i)}$ and $s^{(j)}$}. The resulting matrix is hollow, symmetric, and has positive integers in entries off the diagonal. We input this matrix into Ripser (see {\tt https://live.ripser.org}) and asked for the $H_0$ and $H_1$ barcodes. The result of this experiment can be found in Figure~\ref{fig:poly_barcode}(a).
\end{example}
   \begin{figure}[htbp]
  {\caption{Polyhedra from two sample sizes}\label{fig:poly_sample}}
  {%
    \subfigure[Polyhedra from 20 points]{\label{fig:image-a_sample}%
      \includegraphics[width=0.42\linewidth]{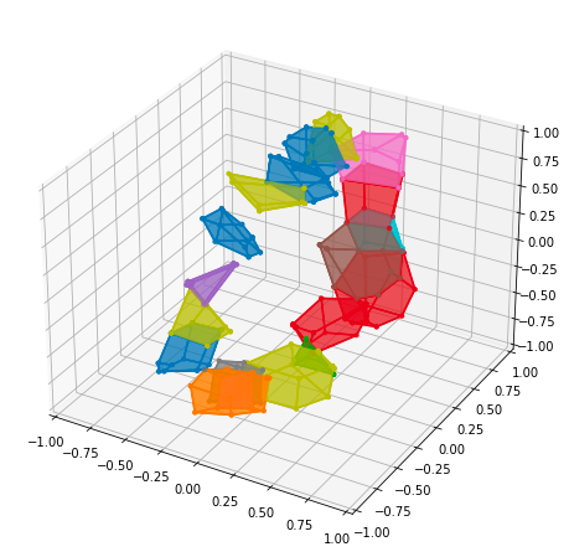}}%
    \qquad
    \subfigure[Polyhedra from 500 points]{\label{fig:image-b_sample}%
      \includegraphics[width=0.42\linewidth]{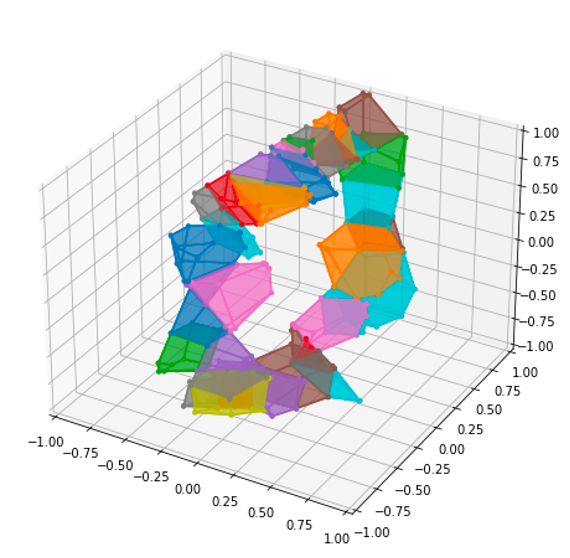}}
  }
\end{figure}
\begin{example} We considered the same circle in $\mathbb R^3$ as in the previous example but sampled at $500$ evenly spaced points $x^{(1)}, x^{(2)}, \dots,$ $x^{(500)}$. A picture of the polyhedra encountered by the $500$ sample points is found in Figure~\ref{fig:poly_sample}(b). We recorded the $500$ binary vectors for each of the polyhedra that was hit by a data point. Only $41$ of the bit vectors were distinct (corresponding to $41$ distinct polyhedra) and we labeled these {$s^{(1)}, s^{(2)}, \dots s^{(41)}$}. We built a $41\times 41$ matrix $F$ by setting $F_{i,j}$ equal to the Hamming distance between {$s^{(i)}$ and $s^{(j)}$}. We input this matrix into Ripser (see {\tt https://live.ripser.org}) and asked for the $H_0$ and $H_1$ barcodes. The result of this experiment can be found in Figure~\ref{fig:poly_barcode}(b).
\end{example}
    \begin{figure}[htbp]
  {\caption{$H_0$ and $H_1$ barcode plots} \label{fig:poly_barcode}}
  {%
    \subfigure[$H_0, H_1$ barcodes for Figure~\ref{fig:poly_sample}(a)]{\label{fig:barcode_a}%
      \includegraphics[width=0.88\linewidth]{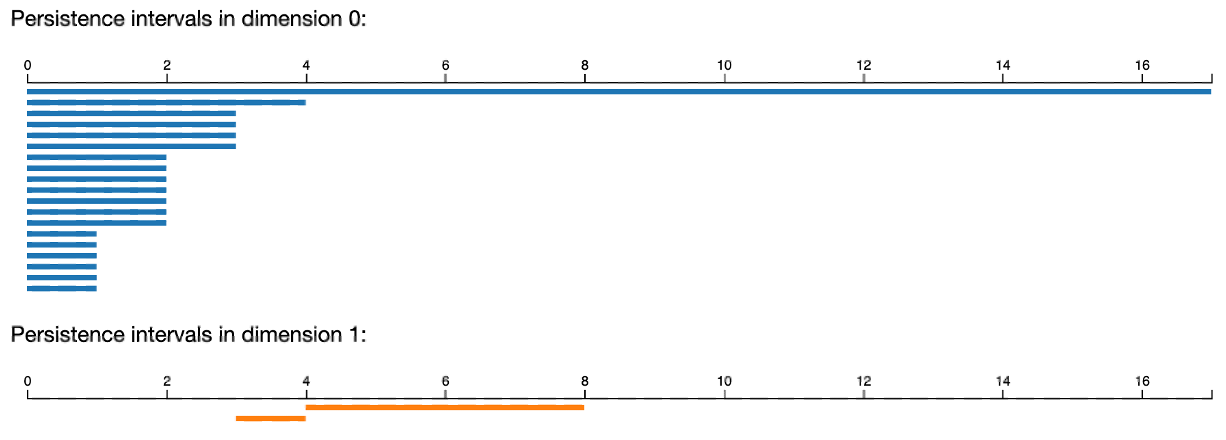}}%
    \qquad
    \subfigure[$H_0, H_1$ barcodes for Figure~\ref{fig:poly_sample}(b)]{\label{fig:barcode_b}%
      \includegraphics[width=0.88\linewidth]{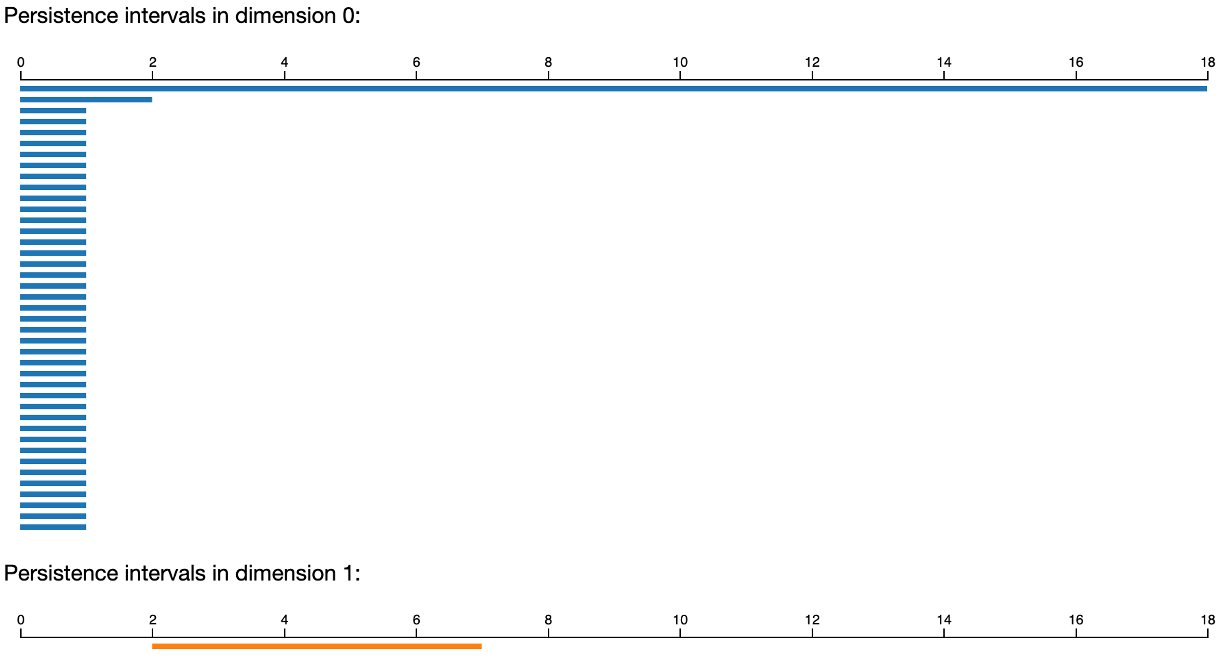}}
  }
\end{figure}
\vspace{-0.1cm}
The $H_0$ barcode from the first example indicates connectivity at Hamming distance 4. The $H_1$ barcode indicates a spurious (i.e. short lived) closed loop occurring at Hamming distance 3. The homological signal of the circle appears (and is quite strong) at Hamming distance 4.
The $H_0$ barcode from the second example indicates connectivity at Hamming distance 2. The $H_1$ barcode indicates a long lived loop beginning at Hamming distance 2.

In the following two examples, we will consider a similar example but utilizing real images from the ImageNet validation dataset (which contains 50K images). We will calculate the Hamming distance between bit vectors of data points via a  much deeper neural network. The network we use is known as
ResNet-50, it is a 50-layer convolutional neural network and was  pre-trained on ImageNet~\cite{deng2009imagenet}. The training images are  224x224x3 thus the input space has dimension 150528. 
It uses ReLU as an activation function on the outputs from the convolutional layers and contains more than 6,000,000 nodes in its many layers.
Like in the previous examples, the activation pattern of each ReLU layer is stored as a bit vector (as defined in 
Section~\ref{subsec:bitvec}).

\begin{example}
\label{example:resnet50_1}
We started with 3 pictures, denoted by $A_1, A_2, A_3$, chosen from the ImageNet validation dataset. They are photos of a miniature poodle, a Persian cat, and a Saluki (see Figure~\ref{fig:images}). Each photo is represented by a 224 x 224 x 3 array. We generated  $50$ data points using the formulas $\sin \theta A_1 + \cos \theta A_2$  and $A_3+\sin \theta A_1 + \cos \theta A_2$, respectively, with $\theta$ consisting of $50$ points uniformly sampled from $[0, 2\pi]$. We calculated the bit vectors via ResNet-$50$ for the $50$ data points and label them as {$s^{(1)}, s^{(2)}, \dots s^{(50)}$}. We built a $50\times 50$ distance matrix $G$ by setting $G_{i,j}$ equal to the Hamming distance between $s^{(i)}$ and $s^{(j)}$. We input this matrix into Ripser which returned the $H_0$ and $H_1$ barcodes. The result of this experiment can be found in Figure~\ref{fig:poly_barcode_resnet50}.
 \end{example}

\begin{figure}[ht]
\centering
  \includegraphics[width=\linewidth]{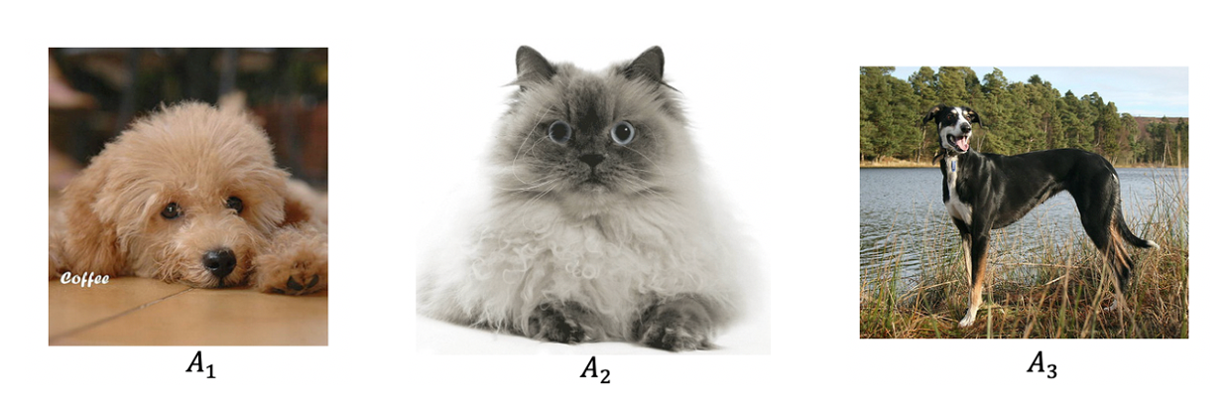}
    \caption{Images of $A_1, A_2,\ \text{and}\  A_3$}
    \label{fig:images}
\end{figure}

We note that there is homological noise (or homological dust) in Figure~\ref{fig:poly_barcode_resnet50}(a) but none in Figure~\ref{fig:poly_barcode_resnet50}(b). In both examples, there is a strong signal representing the circle. The sample points in the second example have all positive entries while the sample points in the first example do not have this property. We were unsure why (or if) this is related to the homological noise but we found it interesting nevertheless. It is a potentially useful feature that the input space to a neural network simultaneously has two kinds of distance measures.  The first derives from the standard Euclidean geometry and the second derives from the coarse geometry implied by the Hamming distance. If one applies an isometry to a data set then its pairwise distance matrix will not change. However, if one applies an isometry with respect to one metric but measure distance via the second metric then one can definitely observe a change. It may be useful to combine information from multiple measurements. One way to carry this out is by using a monotone Boolean function (i.e Boolean functions built using only {\bf and} and {\bf or} operations). The next example illustrates this approach using a small rotation as the isometry in order to produce two distance matrices arising from essentially the same data set.
    \begin{figure}[!ht]
  {\caption{$H_0$ and $H_1$ barcode plots for ResNet-50 \label{fig:poly_barcode_resnet50}}}
  {%
    \subfigure[$H_0, H_1$ barcodes for $\sin \theta A_1 + \cos \theta A_2$]{\label{fig:resnet50_circle_nocenter}%
      \includegraphics[width=0.95\linewidth]{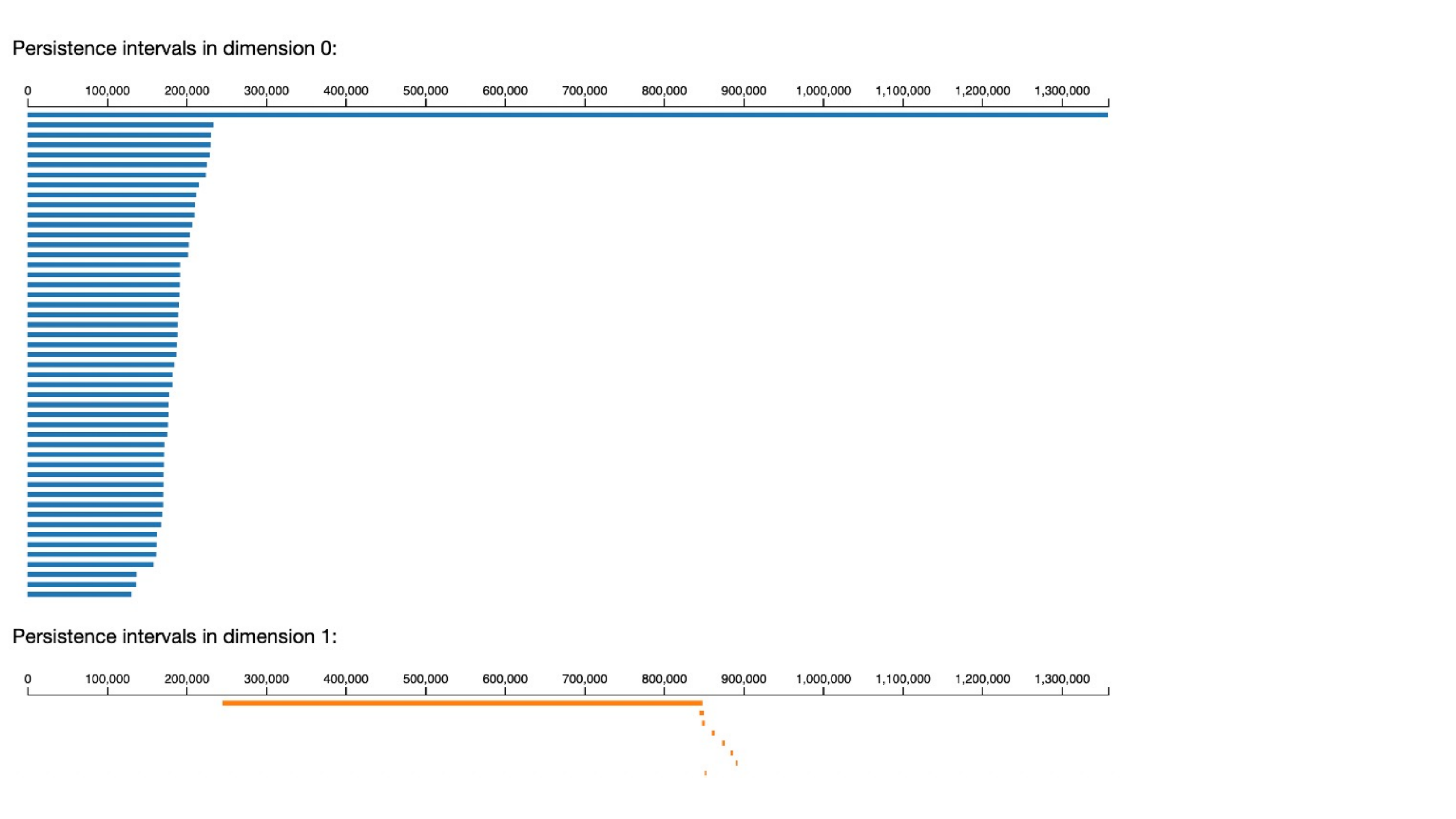}}%
    \qquad
    \subfigure[$H_0, H_1$ barcodes for $A_3+\sin \theta A_1 + \cos \theta A_2$]{\label{fig:resnet50_circle_center}%
      \includegraphics[width=0.95\linewidth]{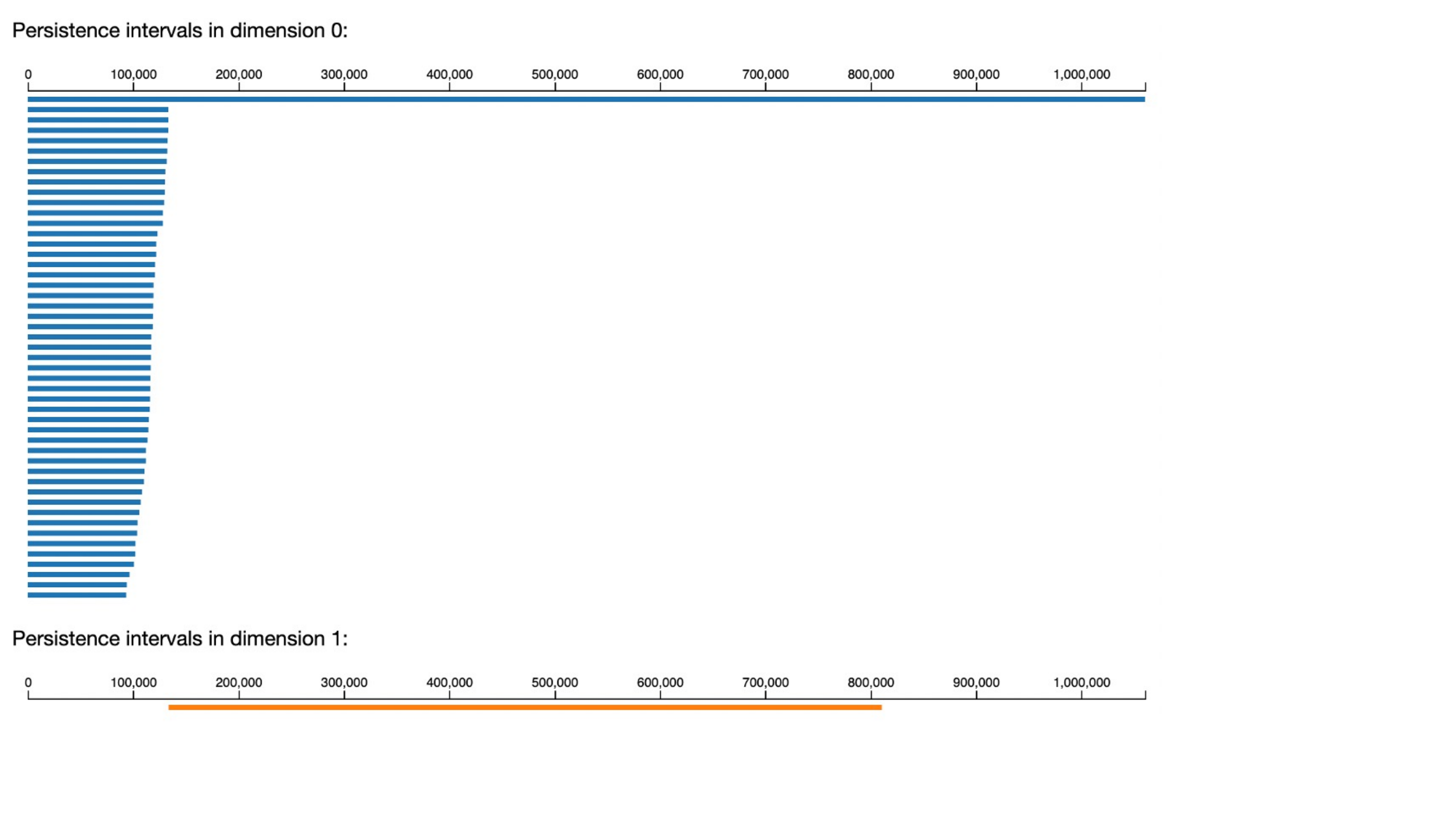}}
  }
\end{figure}
\vspace{-0.1cm}

\begin{example}
\label{example:H1combination}
    Using the same $A_1, A_2$ with Example~\ref{example:resnet50_1}, we uniformly generated 50 points using  $\sin \theta A_1 + \cos \theta A_2$ from $[0, 2\pi]$ and $[1, 2\pi+1]$, respectively. We calculated the bit vectors via ResNet-50 for the $50$ data points sampled from $[0, 2\pi]$  and labelled them as {$s^{(1)}, s^{(2)}, \dots s^{(50)}$}. Similarly, we label the bit vectors of the $50$ data points sampled from $[1, 2\pi+1]$   as {$\bar{s}^{(1)}, \bar{s}^{(2)}, \dots \bar{s}^{(50)}$}. We built four distance matrices $G^1, G^2, G^3,$  and $G^4$ with $G^1_{i,j}$ equal to the Hamming distance between $s^{i}$ and $s^j$, $G^2_{i,j}$ equal to the Hamming distance between $\bar{s}^{i}$ and $\bar{s}^j$, $G^3_{i,j}$ equal to the maximum of the Hamming distance between $s^{i}$ and $s^j$ and between $\bar{s}^{i}$ and $\bar{s}^j$, and $G^4_{i,j}$ equal to minimum of the Hamming distance between $s^{i}$ and $s^j$ and between $\bar{s}^{i}$ and $\bar{s}^j$. The $H_1$ barcodes for the four cases are presented in Figure~\ref{fig:H1barcode}
\end{example}
The distance matrix corresponding to the max function corresponds to requiring that distance matrices $G_1$ and $G_2$ both have their corresponding entry below some threshold. In other words, this corresponds to applying the {\bf and} Boolean function. Similarly, the min function corresponds to applying the {\bf or} Boolean function. Each seems to improve some feature related to the homological noise but it is hard to say which of the two is better. In the next example, we see an overall strengthening of the length of the homological signal for the max function and we see an earlier start of the homological signal for the min function. Homological noise did not make an appearance in the next example.
\begin{example}
     Using the same $A_1, A_2, A_3$ with Example~\ref{example:resnet50_1},  we uniformly generated 50 points using  $A_3+\sin \theta A_1 + \cos \theta A_2$ from $[0, 2\pi]$ and $[1, 2\pi+1]$, respectively. As in Example~\ref{example:H1combination}, we generated 4 different distance matrices $\hat{G}^1$, $\hat{G}^2$, $\hat{G}^3$, and $\hat{G}^4$, respectively. The $H_1$ barcodes for the four cases are presented in Figure~\ref{fig:H1barcode_shift}.
\end{example}
\begin{figure}[ht]
\centering
  \includegraphics[width=0.9\linewidth]{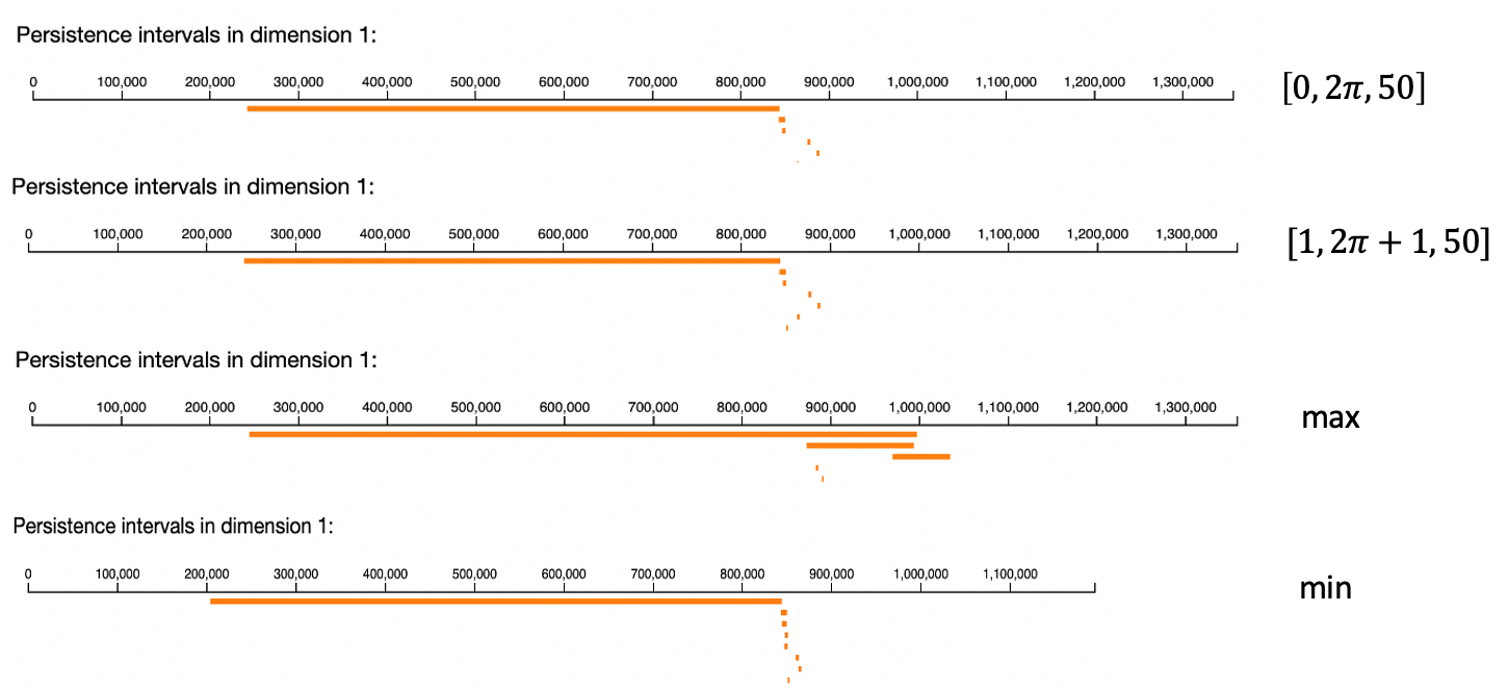}
    \caption{$H_1$ barcodes for $G^1$, $G^2$, $G^3$, and $G^4$}
    \label{fig:H1barcode}
\end{figure}
\begin{figure}[ht]
\centering
  \includegraphics[width=0.9\linewidth]{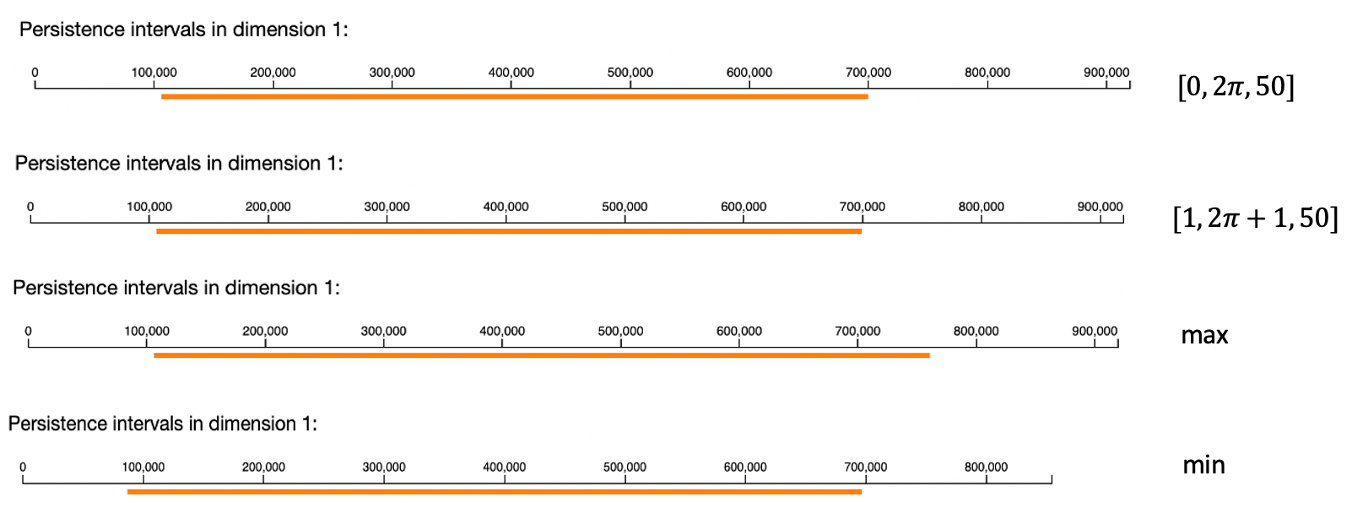}
    \caption{$H_1$ barcodes for $\hat{G}^1$, $\hat{G}^2$, $\hat{G}^3$, and $\hat{G}^4$}
    \label{fig:H1barcode_shift}
\end{figure}

\section{Conclusion and Future Work}\label{sec:concl} 
A ReLU feedforward neural network induces a finite polyhedral decomposition of the input space. The corresponding dual graph represents this decomposition, where vertices correspond to polyhedra and edges represent shared facets.  Data in the input space gets mapped to vertices of this graph. The geometry of the graph can be exploited, via persistent homology, to detect homological signals of manifolds from points sampled from the manifold. Many techniques in data analysis build on the premise that data, sharing a collection of common features, tends to be well approximated by a low dimensional geometric object. It is conjectured that the coarseness of the polyhedral decomposition can be helpful in dealing with the noise present in many data sets (maybe by including monotone Boolean functions to combine multiple distance matrices). In future work, we hope to extend the examples in this paper with the goal of detecting manifolds, known to exist in data, such as $SO(3)$, higher dimensional tori, and various fiber bundles involving these manifolds. A torus example can be found in Appendix~\ref{appendix:example}.

Neural networks of the kind considered in this paper have a finite number of ``polyhedral resources". The training function, the distribution of training data, the network architecture, and the training method are some of the factors in how the network utilizes its polyhedra. In the networks we considered, we observed the general pattern that there was a dense collection of smaller volume polytopes near the training data, larger volume polytopes in the general vicinity but away from the training data, with a shell of unbounded polyhedra surrounding the polytopes. You can see hints of this tendency in Figure \ref{fig:visualization}(a). In future work, we hope to make more precise characterizations in the direction of this observation. In particular, we are eager to extend this work to depictions that more accurately quantify the number of polytopes that emerge from a neural network and their distributions of volumes, the ``polytopes landscape''. 
\section{Acknowledgements}
This work is partially supported by the United States Air Force under Contract No. FA865020C1121 and the DARPA Geometries of Learning Program under Award No. HR00112290074.

\bibliographystyle{icml2023}
\bibliography{sample-base, pmlr-sample_revised}

\newpage
\appendix
\onecolumn

\section{Algorithms}
\label{appendix:algorithms}
\begin{algorithm} 
\caption{Brute-Force Method} 
\label{alg:bruteforce}
\begin{algorithmic} 
    \REQUIRE A pretrained $(L+1)$-layer ReLU FFNN with weights and biases $\{W_i, b_i\}_{i=1}^L$ for all $L$ hidden layers and the number of nodes in all hidden layers equal to $h$
    \ENSURE $\mathcal{P}$, the set of bit vectors that determine all polyhedra in the input space, and the size of $\mathcal{P}$, the number of polyhedra {in the input space}.  
    \STATE 1) Initialize an empty set $\mathcal{P}$ to store all bit vectors that determine a polyhedron.
\STATE 2)  For $j=1$, generate $A^{(j)}x\leq c^{(j)}$ based on $s^{(j)}$ and $\{W_i, b_i\}_{i=1}^L$ using (\ref{ineq:polytopecond}).
\STATE 3) Solve the linear program:
        \begin{equation}
        \label{eq:LP}
\begin{aligned}
\text{minimize } \quad &0 \\
\text{subject to }\quad &A^{(j)}x\leq c^{(j)}.
\end{aligned}   
\end{equation} 
If the optimal value is $0$ (as opposed to `infeasible' if using \texttt{scipy.optimize.linprog()} for example), then {there is at least one feasible solution to (\ref{eq:LP}), so} the bit vector $s^{(j)}$ determines a polyhedron, $p^{(j)}$. If $p^{(j)}$ is full-dimensional, add $s^{(j)}$ to $\mathcal{P}$.
\STATE 4) Repeat 2) and 3) for  $j=2, \ldots, 2^h$.
\end{algorithmic}
\end{algorithm}
\begin{algorithm} 
\caption{Traversal-and-Pruning Method} 
\label{alg:traversal_pruning}
\begin{algorithmic} 
    \REQUIRE A pretrained $(L+1)$-layer ReLU FFNN with weights and biases $\{W_i, b_i\}_{i=1}^L$ for all $L$ hidden layers and a random point $x^{(0)}$ in the input space $\mathbb{R}^m$. Denote by $P(x^{(0)})$ the polytope to which $x^{(0)}$ belongs and by $s(x^{(0)})$ the specific bit string that uniquely identifies $P(x^{(0)})$.
    \ENSURE $\mathcal{P}$, the set of bit vectors that determine all polyhedra in the input space, and the size of $\mathcal{P}$, the number of polyhedra {in the input space}.  
    \STATE 1) Initialize an empty set $\mathcal{P}$ to store all bit vectors that determine a polyhedron and define a label for each stored bit vector with 1 denoting if we have found all polyhedra adjacent to $P(x^{(0)})$ and 0 otherwise.
\STATE 2)  Calculate the bit vector $s(x^{(0)})$ using (\ref{eq:outputderi}) and (\ref{eq:bitvec}), add $s(x^{(0)})$ to the set $\mathcal{P}$, and {give a label of 0 to $s(x^{(0)})$}.
\STATE 3) Traverse all the polyhedra adjacent to $P(x^{(0)})$, save {their bit vectors} to $\mathcal{P}$, and update {these bit vectors'} labels using the following process:
  \begin{itemize}
      \item Find the active bits of $s(x^{(0)})$ by solving the LP model (13) in \cite{vincent_reachable_2021}. {Without loss of generality, say there are $q$ active bits.}
      \item {Flip each of these $q$ active bits of $s(x^{(0)})$ (one at a time), producing $q$ new bit strings $\{\hat{s}^{(j)} : 1 \leq j \leq q \} \eqqcolon \hat{\mathcal{S}}$, which all differ from $s(x^{(0)})$ by only one bit. Add the elements of $\hat{\mathcal{S}}$ to $\mathcal{P}$ that are not already in $\mathcal{P}$.}
      
      \item Label the added bit vectors with 0 and {switch the label of $s(x^{(0)})$ from 0 to 1, indicating that we have found all of the neighbors of $P(x^{(0)})$ and have added their bit vectors to $\mathcal{P}$.} 
  \end{itemize}
\STATE 4) Repeat 3 for bit vectors in $\mathcal{P}$ with label 0 until all the bit vectors in $\mathcal{P}$ have label $1$.
\end{algorithmic}
\end{algorithm}

\section{2-D Polyhedral Decomposition Up Close}\label{appendix:2Dpoly}
With Figure \ref{fig:decomp_with_dicts} and relevant discussion, we demonstrate qualities of the bit vectors that label the distinct regions of the polyhedral decomposition derived from model (\ref{example:model1}). 
We can see that the binary vectors produced by model (\ref{example:model1}) consist of 6 elements or bits, 3 from the first layer of the model and another 3 from the second layer. The collection of the lines of the polyhedral decomposition are defined by the weight matrices and bias vectors of the model.   In the bdictionary values 

It can be verified that neighboring polyhedra (polyhedra that share a 1-dimensional face) only differ in one bit. More specifically, they differ in the entry associated with the inequality that defines their shared edge. For example, looking at the top dictionary in Figure \ref{fig:decomp_with_dicts}, we see that the binary vectors from polytopes 22 and 23 differ in their first bit. The first bit in the binary vector in the bottom dictionary is red, so the edge separating polytopes 22 and 23 is red. 

Additionally, we clarify the claim that a polyhedral decomposition's dual graph (with vertices corresponding to the polyhedra and edges connecting vertices if their corresponding polyhedra share an edge) is bipartite. In our 2-dimensional example, we have two vertex sets $\mathcal{V}_1$ and $\mathcal{V}_2$ where are $\mathcal{V}_1 =\{2,4,6,7,8,9,10,12,15,17,19,22,24,26, \text{ and the unlabeled center triangle} \}$ and $\mathcal{V}_2 = \{ 1,3,5,11,13,14,16,18,20,21,23, \text{and } 25 \}$ because no two polyhedra in the same set are neighbors.

\begin{figure}[h]
\centering
    \includegraphics[width = \linewidth] 
    {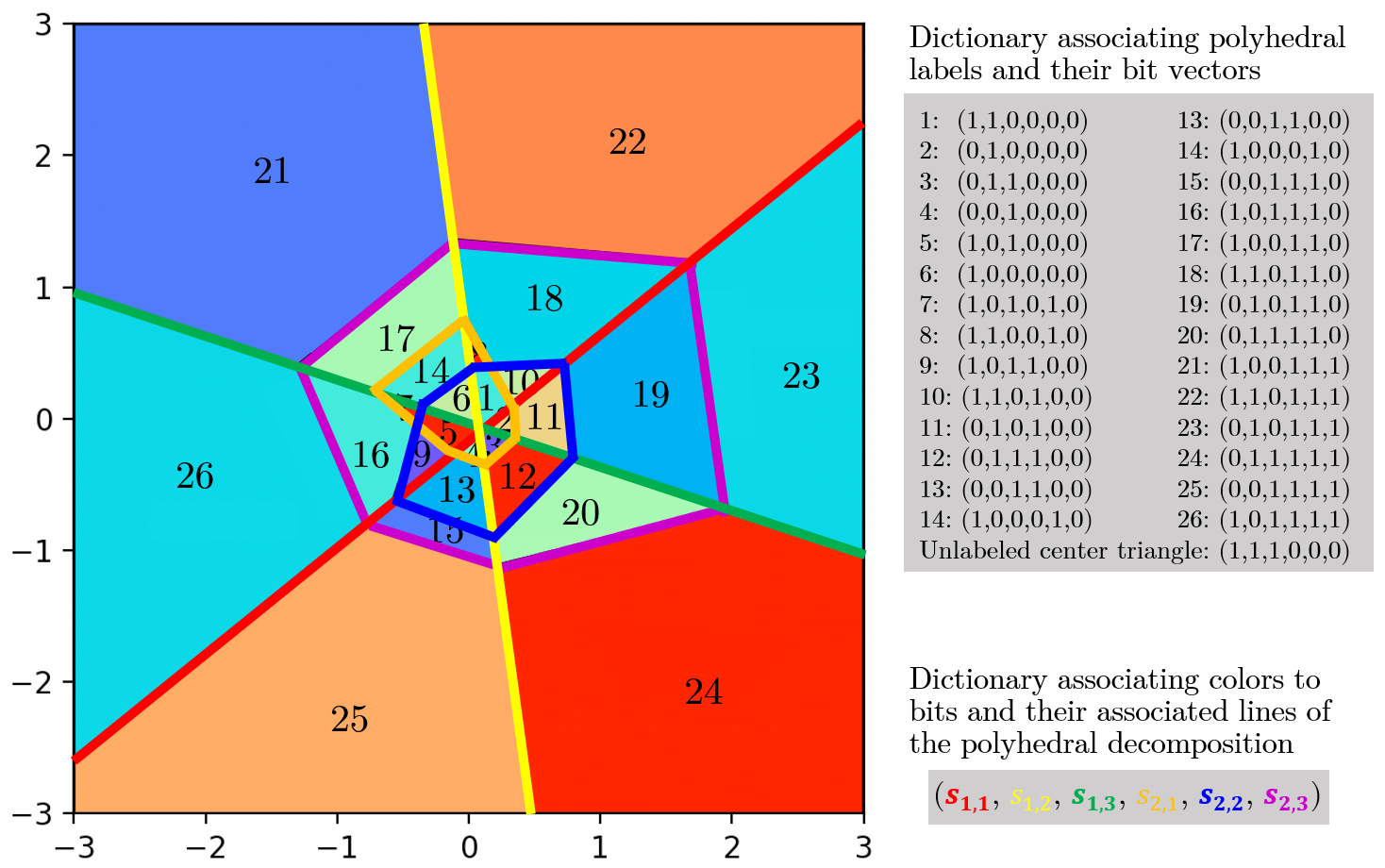}
    \caption{The 2-dimensional polyhedral decomposition with a look-up table giving the binary vectors for each polytope and color-coding that associates each bit to lines in the polyhedral decomposition using the notation introduced in (\ref{eq:bitvector}). }
    \label{fig:decomp_with_dicts}
\end{figure}

\section{Example}
\label{appendix:example}
\begin{example}
In this example, we started with 5 randomly generated orthogonal tensors of size $224 \times 224 \times 3$ (denoted $A_1, \dots, A_5$).
 We made $100$ data points generated using the formula $A_5+\alpha*(\sin \theta_1 A_1 + \cos \theta_1 A_2 +\sin \theta_2 A_3 + \cos \theta_2 A_4)$, for $\theta_1, \theta_2$ being $100$ points uniformly sampled from $[0, 2\pi]\times [0, 2\pi]$. 
 We input this matrix into Ripser and asked for the $H_1$ and $H_2$ barcodes. 
 The result of this experiment can be found in Figure~\ref{fig:resnet50_torus_generated}.
\end{example}
\begin{figure}[ht]
\centering
     \includegraphics[width = 1\linewidth]
     {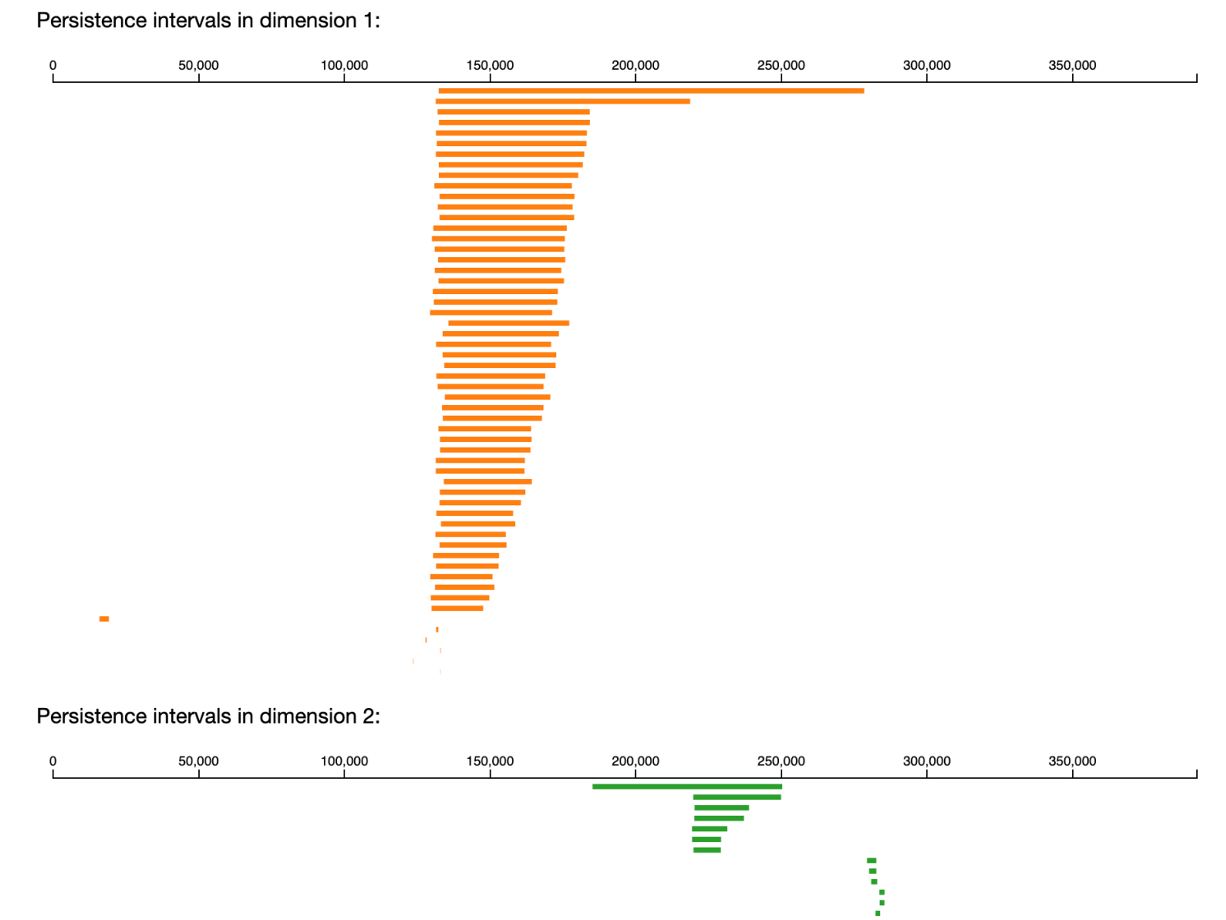}
    \caption{$H_1$ and $H_2$ barcode plots of randomly generated orthogonal tensors.}
    \label{fig:resnet50_torus_generated}
\end{figure}
\end{document}